\documentclass{amsart}
\usepackage[utf8]{inputenc}
\usepackage{amsmath,amssymb,bm,amsthm,xcolor}

\newcommand{\bb}[1]{\bm {#1}}
\newcommand{\beq}{\begin{equation}}
\newcommand{\eeq}{\end{equation}}
\newcommand{\dfdx}[2]{\frac{\partial {#1}}{\partial{#2}}}
\newcommand{\dfdxy}[3]{\frac{\partial^2 {#1}}{\partial{#2} \partial{#3}}}
\newcommand{\R}{\mathbb{R}}
\newcommand{\pref}[1]{(\ref{#1})}

\newtheorem{theorem}{Theorem}

\begin{document}

\title{When is the inverse  of an invertible convex function itself convex?}
\author{Robert Planqu\'e}
\address{
Department of Mathematics, Vrije Universiteit Amsterdam\\
De Boelelaan 1111, 1081 HV, Amsterdam, the Netherlands}
\email{r.planque@vu.nl}
\date{}

\subjclass{26B25}
\dedicatory{This work was the direct result of the long-standing collaboration with Joost Hulshof and Frank Bruggeman. In particular, the author wishes to thank Hulshof for his unrelenting efforts, including providing the general view on concavity and convexity in this note. Andr\'e Ran encouraged me to publish this result as a separate note, and Lionel Thibault provided very constructive feedback that considerably improved the manuscript.}
\keywords{invertible convex function, strongly convex function}

\begin{abstract}
  We provide a sufficient condition for an invertible (locally strongly) convex vector-valued  function on 
$\mathbb{R}^N$ to have a (locally strongly) convex inverse. We show under suitable conditions that if the gradient of each component of the inverse has negative entries, then this inverse is (locally strongly) convex if the original is.  \end{abstract}

\maketitle

In this short note, we provide a sufficient condition for an invertible convex function on $\R^N$ to have a convex inverse. To the best of our knowledge, it has not appeared before in the literature, and since this type of result is often buried within longer articles, we thought it would be more easily found when isolated.  

Let $D$ be a non-empty open convex set in $\mathbb{R}^N$ and $f: D \to \mathbb{R}$ be a real-valued function on $D$. The function $f$ is called {\em strictly convex} on a convex subset $C$ of $D$ in the literature if 
$$
  f((1-t){\bm{u}}+ t{\bm{v}}) < (1-t)f({\bm{u}}) + tf({\bm{v}}) \quad 
	\text{for all}\; t\in (0,1)\;\text{and}\;{\bm{u}},{\bm{v}}\in C \;\text{with}\;
	{\bm{u}}\neq {\bm{v}}.
$$
An important subclass is that of strongly convex functions. 
We recall that $f$ is {\em strongly convex} on a convex subset $C$ of $D$ (see, e.g., \cite[p. 565]{Rock-Wets}, \cite[p. 430]{Thib}) whenever there exists a real $\sigma >0$ such that for any 
${\bm{u}},{\bm{v}} \in C$ and any $t\in (0,1)$ 
\begin{equation}\label{eq-def-StConv}
    f((1-t){\bm{u}}+ t{\bm{v}})\leq (1-t)f({\bm{u}}) + tf({\bm{v}}) 
- \frac{1}{2}\sigma{t(1-t)}\|{\bm{u}}-{\bm{v}}\|^2, 
\end{equation}	
where $\|\cdot\|$ denotes the canonical Euclidean norm of $\mathbb{R}^N$. By \cite[Exercise 12.59(c)]{Rock-Wets} and  \cite[Proposition 3.273]{Thib} the latter inequality for strong convexity of $f$ on $C$ is equivalent to the 
mere convexity on $C$ of the function $f -\frac{1}{2}\sigma\|\cdot\|^2$. The strong convexity evidently implies the strict convexity, but the converse is known to fail (see, e.g., \cite[p. 430]{Thib}). The function $f$ is {\em locally strongly convex} on $D$ if for each ${\bm{x}}$ in $D$ there is an open convex neighborhood $U \subset D$ of ${\bm{x}}$ on which $f$ is strongly convex.

It is worth noticing that the real-valued function $f$ is convex on $D$ whenever it is locally strongly convex on $D$ since it is known that any locally convex real-valued function on a convex set is convex on this set.  
We declare a vector-valued function  $f: D\to \mathbb{R}^N$ to be {\em locally strongly convex} on $D$ if each of its components is. If $N=1$ and $C$ is an open interval $I$ over which the real-valued function $f$ is differentiable, then the strict convexity of $f$ on $I$ (see, e.g., \cite[Theorem 2.13]{Rock-Wets}, \cite[p. 98]{Thib}) 
 is equivalent to the increasing property on $I$ of the derivative 
$f'$, that is, $f'(u) < f'(v)$ for $u< v$ in $I$. Regarding the above inequality of strong convexity, by \cite[Exercice 12.59(c)]{Rock-Wets} and  \cite[Proposition 3.273]{Thib}, it means under the differentiability (resp. twice differentiability) assumption of $f$ on $I$ that 
$f'(v) \geq f'(u) +\sigma(v-u)$ for $u \leq v$ in $I$ (resp. $f''(u) \geq \sigma$ 
for $u\in I$). Then, for $f$ of class $C^2$ on the open interval $I$, the local strong convexity of 
$f$ on $I$ is equivalent to the positivity on $I$ of the second derivative $f''$, that is, 
$f''(u) > 0$ for $u\in I$.

For a scalar  $C^2$ locally strongly convex function $f:\mathbb{R}\to \mathbb{R}$ with $f'(x)\neq 0$ that is invertible (perhaps when restricted to some open interval), it is easy to see that the inverse is convex if $g'(y) < 0$. Let $g$ be the inverse of $f$, and denote $f(x)=y$ so that $g(y)=x$.

Assume in accordance with the above that $f'(x) \ne 0$, $f''(x)>0$, $g'(y)<0$. Since $g(f(x))=x$, implicit differentiation gives
\[
g'(f(x)) f'(x) = 1.
\]
Differentiating once more yields
\[
g''(f(x)) f'(x) ^2 + g'(f(x)) f''(x)=0,
\]
so that
\[
g''(f(x)) = - \frac{g'(f(x)) f''(x)}{f'(x)^2} > 0.
\]
A simple example is $f(x)=1/x$ for $x>0$, which is its own inverse, and is locally strongly convex. A second example is $f(x) = e^{-x}$, with inverse $-\ln x$, both of which are locally strongly convex. Not every convex function has a convex inverse, as the example $f(x)=e^x$ shows.

To generalise this result to higher dimensions, we need a few basic facts. 

 Let $D$ be an open convex set of $\mathbb{R}^N$ and $f:D \subset \mathbb{R}^N \to \mathbb{R}$ be a real-valued function.  By what is recalled above, $f$ is locally strongly convex on $D$ if and only if for each ${\bm{x}}\in D$ there is an open convex neighborhood $U$ of $\bm{x}$ and a real 
$\sigma_{\bm{x}} >0$ such that the function 
$f -\frac{1}{2}\sigma_{\bm{x}}\|\cdot\|^2$ is convex on $U$. Assuming that $f$ is of class $C^2$, the function $f$ is then strongly locally convex on $D$ if and only if its Hessian matrix is positive definite 
at each point of the open convex set $D$.

Recall that two $n\times n$ matrices $A$ and $B$ are congruent if there exists an invertible matrix $W$ such that $W^T A W = B$ \cite{hornjohnson.85}. Congruence is an equivalence relation. If $A$ is positive definite, then all matrices that are congruent to $A$ are also positive definite \cite{hornjohnson.85}:
\[
\bb x^T (W^T A W) \bb x = (W \bb x)^T A (W \bb x) >0.
\]
Moreover, 
the set of positive definite matrices forms a positive cone, i.e., positive linear combinations of positive definite matrices are again positive definite \cite[p. 399]{hornjohnson.85}.


In the theorem below we show that under some differentiability assumption, also in higher dimensions, the condition that the component functions of the inverse have negative gradients suffices for them to be locally strongly convex. Since each component is locally strongly convex, so is the full inverse function.

Below 
$\mathbb{R}^N_{< 0}$ stands for $\mathbb{R}^N_{< 0}:=(-\infty,0)^N$ and $\delta_{mn}$ is defined by 
$$ \delta_{mn}=\left\{ 
\begin{array}{l} 
1 \quad\text{if}\; m=n \\
0 \quad\text{otherwise.}
\end{array}
\right. 
$$

\begin{theorem}
\label{thm:inverseconvex}
Let $D$ and $R$ be non-empty open convex sets of $\R^N$. Let $f : D \subset \R^N \to R \subset \R^N$ be $C^2$, locally strongly convex and invertible, and let $g : R \to D$ be its inverse. Let $\bb x \in D$ and $\bb y \in R$ be such that $f(\bb x)=\bb y$, $g(\bb y)=\bb x$. If for each component $g_m$, $\nabla g_m(\bb y) \in \R^N_{<0}$, then $g$ is locally  strongly convex on $R$.
\end{theorem}

\begin{proof}
 Since $g(f(\bb x)) = \bb x$ for each $\bb x \in D$, differentiation of the $m$-th component of $g$ with respect to $x_n$ gives
 \beq
 \label{eq:1}
\sum_{i=1}^N \dfdx {g_m}{y_i}(f(\bb x)) \dfdx {f_i}{x_n}(\bb x) = \delta _{mn}.
\eeq
This can also be written as
\[
\nabla g_m(\bb y)^T \dfdx {f}{x_n} = \delta_{mn}
\]
for each $m,n=1,\ldots,N$ and in which $\dfdx{f}{x_n}$ denotes the vector containing derivatives of all component functions of $f$ with respect to $x_n$. 
 For the entire function $g$, these results may be summarized by 
\beq
\label{eq:DfDg}
Dg(\bb y) Df(\bb x) = I_N,
\eeq
showing that the Jacobi matrices of $f$ and $g$ are each others inverses.
Focusing again on one component $g_m$, further differentiation of \pref{eq:1} with respect to $x_j$ yields
\begin{multline}
\sum_{l=1}^N \dfdxy {g_m}{y_1}{y_l} \dfdx {f_l}{x_n} \dfdx {f_1}{x_j} + \dfdx {g_m}{y_1} \dfdxy {f_1} {x_j}{x_n} + \cdots +\\
\sum_{l=1}^N \dfdxy {g_m}{y_N}{y_l} \dfdx {f_l}{x_n} \dfdx {f_N}{x_j} + \dfdx {g_m}{y_N} \dfdxy {f_N} {x_j}{x_n} = 0,
\end{multline}
where we have suppressed the evaluation of all these functions at $\bb x$ or $\bb y$ in the appropriate places. We may abbreviate this to
\beq
\label{eq:HgHf}
H_{g_m}\left(\dfdx f{x_j},\dfdx f{x_n}\right) + \nabla g_m^T \dfdxy{f}{x_j}{x_n} = 0,\quad j,n=1,\ldots,N.
\eeq
On the left we see the Hessian of $g_m$ applied as a bilinear form to the directions $\dfdx f{x_j}$ and $\dfdx f{x_n}$. These directions $\dfdx f{x_j}$ form a basis since collectively they form the Jacobian $Df$, which is invertible. 
We can also write \pref{eq:HgHf} as
\beq
\label{eq:DfHDf}
Df^T H_{g_m} Df 
= -\left(\sum_{k=1}^N \dfdx {g_m}{y_k} \dfdxy {f_k}{x_j}{x_n}\right)_{j,n=1\ldots,N}.
\eeq
The right hand side is a linear combination of Hessians of each component, $H_{f_k}$, with coefficients from $\nabla g_m$. Since $f$ is locally strongly convex, so are all components $f_k$. The Hessian of such a component $H_{f_k}$ is thus positive definite, and if we assume that $\nabla g_m$ is a vector with negative entries, then the right hand side of \pref{eq:DfHDf} is again positive definite. But then \pref{eq:DfHDf} states that $H_{g_m}$ is congruent to a positive definite matrix, making $H_{g_m}$ itself positive definite.
\end{proof}

It is clear that the conditions in the theorem can be relaxed somewhat, for instance by requiring that the component functions $g_m$ only decrease in some variables, and do not vary in others. The conditions are  sufficient but are far from necessary. For example, the identity function on $\R^n$ is convex (but not strict convex); its inverse, which is again the identity, is nowhere decreasing, however.

\section*{Origin of the problem}

Our interest in the question whether inverses of convex vector-valued functions are themselves convex was raised in the following context. More detail and motivation may be found in \cite{planque.18,planque.23}.

Consider a biochemical pathway
\beq
\label{linchain}
\underline{x_0} \stackrel{e_1}\rightleftharpoons x_1
\stackrel{e_2}\rightleftharpoons x_2
\stackrel{e_3}\to\emptyset,
\eeq
in which so-called metabolites $X_1$ and $X_2$, with concentrations $x_1$ and $x_2$ are converted from a metabolite $X_0$ with fixed concentration $\underline{x_0}$ by enzymes $E_1$, $E_2$ and $E_3$, with corresponding concentrations $e_1$, $e_2$ and $e_3$. The emptyset at the end of the pathway indicates that the product made by the last enzyme does not influence the reaction rate. Clearly, all concentrations are positive.

Treating the enzyme concentrations as parameters, a typical dynamical system for the evolution of the metabolite concentrations $x_1$ and $x_2$ is given by
\beq
\label{dynsys}
\begin{aligned}
    \dot x_1 &= e_1 f_1(\underline{x_0},x_1) - e_2 f_2(x_1,x_2),\\
    \dot x_2 &=  e_2 f_2(x_1,x_2) - e_3f_3(x_2).
\end{aligned}
\eeq
The functions $f_1,\ldots,f_3$  describe how the reaction rate of the corresponding catalytic enzyme is influenced by the concentrations of its substrates and products. In this case, a typical choice is (reversible) Michaelis-Menten kinetics,
\beq 
\label{MMold}
f_i = k_{cat,i}\frac{x_{i-1} - \frac{x_i}{K_{eq,i}}}{\frac{x_{i-1}}{K_{i,i-1}} + \frac{x_i}{K_{i,i}} + 1},\quad i=1,2,3,
\eeq
with $x_0 = \underline{x_0}$ and $x_3=0$. All parameters here are strictly positive. For this illustration, we rescale variables and parameters to simplify \pref{MMold} to 
\beq
\label{MM}
f_i(x_{i-1},x_i) = \frac{x_{i-1} - x_i}{a_i x_{i-1} + b_i x_i + c_i}.
\eeq
The dynamical system \pref{dynsys} with \pref{MM} is globally stable on the state space of positive metabolite concentrations $\R^2_{\ge 0}$ \cite{smillie.84}. If $\underline{x_0}>0$, any steady state solution satisfies
\beq
\label{stst}
e_1 f_1(\underline{x_0},x_1)
= e_2 f_2(x_1,x_2) = e_3f_3(x_2) = J,
\eeq
with 
\beq
\label{x012}
\underline{x_0} > x_1 > x_2 > 0
\eeq
and $J>0$ is the so-called flux through the pathway. For any positive set of enzyme concentrations $\bb e$ and constant $\underline{x_0}$, we thus find a flux $J=J(\bb e)$ by solving \pref{stst} on \pref{x012}. 

From \pref{stst}, it follows that $J$ is 1-homogeneous in $\bb e=(e_1,e_2,e_3)$: $J(\lambda \bb e) = \lambda J(\bb e)$ for any positive $\lambda$. The quantity
\[
\frac{J}{e_1+e_2+e_3},
\]
is thus 0-homogeneous. It is the flux per unit total enzyme expended to sustain this flux, and called the specific flux. Optimizing this quantity is central to the study the metabolic behaviour of single-celled organisms such as bacteria and yeasts \cite{bruggeman.20}. 

For the simple linear pathway \pref{linchain} we are therefore interesting in studying
\beq
\label{maxprob}
\begin{aligned}
\max_{x_1,x_2,e_1,e_2,e_3}\Big\{ \frac{J}{e_T}\ |\ J & = e_1f_1(\underline{x_0},x_1) = e_2f_2(x_1,x_2) = e_3f_3(x_2),\\ 
&\quad\quad\quad e_i\ge 0,\ i=1,2,3,\\
&\quad\quad\quad  e_1+e_2+e_3=e_T,\\ 
&\quad\quad\quad  \underline{x_0} > x_1 > x_2 > 0  \Big\}.
\end{aligned}
\eeq
Of course, instead of maximising $J/e_T$, we may minimize $e_T/J$. This actually simplifies \pref{maxprob}: since $e_i = J/f_i(\bb x)$  in steady stat, this is equivalent to minimising 
\beq
\label{obj}
\frac{1}{f_1(\underline{x_0},x_1)} + \frac 1{f_2(x_1,x_2)} + \frac{1}{f_3(x_2)},
\eeq
and \pref{maxprob} reduces to a problem in which only $x_1$ and $x_2$ enter as optimisation variables,
\beq
\min_{x_1,x_2}\Big\{ \frac{1}{f_1(\underline{x_0},x_1)} + \frac 1{f_2(x_1,x_2)} + \frac{1}{f_3(x_2)}\ |  
\  \underline{x_0} > x_1 > x_2 > 0  \Big\}.
\eeq
With the choice \pref{MM} the function \pref{obj} is locally strongly convex in logarithmic variables (i.e., setting $y_i = \log x_i$) on the induced domain $\underline{y_0} > y_1 > y_2 > -\infty$ \cite{noor.16,planque.18}. 

From a cellular perspective, however, it is more relevant to understand how the enzyme concentrations impact \pref{maxprob}. Single-celled organisms control their lives by making these enzymes through gene expression; the metabolites involved in the metabolic pathway they employ are not under their direct control. To optimize the flux through the pathway, ideally the specific $J(\bb e)/e_T$ would need to be concave. For the particular example of the small linear chain of reactions, this may in fact be proven directly, by implicit differentiation of the steady state equations \pref{stst} to second order with respect to the enzyme concentrations, and studying the determinant of the Hessian of $J$. For more complicated pathways, however, this method is fruitless.

To provide a more general technique to link convexity of \pref{obj} in logarithmic variables to concavity of $J(\bb e)/e_T$, we first introduce 
\[
z_i = \frac{e_i}J,\quad i=1,2,3,
\]
so that the steady state equations \pref{stst} read
\beq
\label{zstst}
z_1 f_1(\underline{x_0},x_1)
= z_2 f_2(x_1,x_2) = z_3f_3(x_2) = 1.
\eeq
In logarithmic variables $\bb y$ for $\bb x$, $\bb z = (z_1,z_2,z_3)$ and $(\underline{y_0},y_1,y_2)$ are thus related by
\beq
z_1 = \frac 1{f_1(e^{\underline{y_0}},e^{y_1})},\
z_2 = \frac 1{f_2(e^{y_1},e^{y_2})},\
z_3 = \frac 1{f_3(e^{y_2})}.
\eeq
With the choice \pref{MM}, these equations may be explicitly inverted,
\beq
\label{inverses}
y_2 = \log F_3(z_3),\quad y_1 = \log F_2(z_2,z_3),\quad
\underline{y_0} = \log F_1(z_1,z_2,z_3),
\eeq
with the $F_i$ defined by
\begin{align*}
F_3(z_3) &= \gamma_3 M_3(z_3),\\
F_2(z_2,z_3) &= \gamma_2 M_2(z_3) + \gamma_3 M_3(z_3)(1+ M_2(z_2)) \\
F_1(z_1,z_2,z_3) &= \gamma_1 M_1(z_1) +
\gamma_2 M_2(z_2) + \gamma_3 M_3(z_3(1+M_2(z_2))(1+M_1(z_1)),
\end{align*}
and
\[
\gamma_i = \frac{c_i}{z_i - a_i},\quad
M_i = \frac{a_i+b_i}{z_i-a_i},\quad i=1,2,3.
\]
The natural domain of these $F_i$ functions (induced by \pref{x012}) is $z_i > a_i$, $i=1,2,3$. On this domain, it may readily be checked that the inverse functions in \pref{inverses} have negative gradients. Theorem 1 thus shows that  
they are locally strongly convex functions. 
From this, it may also be concluded that the $F_i$ themselves are locally strongly convex on $z_i > a_i$.


By studying the relation between the Hessians of $F_1(\bb z)$ and  $J(\bb e)$, one can now show that $J(\bb e)/e_T$ is concave for all pathways for which the original objective function \pref{obj} is locally strongly convex in logarithmic variables \cite{planque.23}. The condition that the inverse functions have negative gradients turns out to be natural in the more general pathway setting. This technique thus works for many complicated pathways, even when the inverse functions $\log F(\bb z)$ cannot be explicitly computed.


\begin{thebibliography}{1}
\providecommand{\natexlab}[1]{#1}
\providecommand{\url}[1]{\texttt{#1}}
\expandafter\ifx\csname urlstyle\endcsname\relax
  \providecommand{\doi}[1]{doi: #1}\else
  \providecommand{\doi}{doi: \begingroup \urlstyle{rm}\Url}\fi

\bibitem[1]{Rock-Wets}  R.~T. Rockafellar and R. J.-B. Wets. {\em Variational Analysis},
Springer-Verlag, Berlin, 1998.  
		
\bibitem[2]{Thib} L. Thibault. {\em Unilateral Variational Analysis in Banach Spaces. Part I: General Theory}, 
World Scientific, 2023. 


\bibitem[3]{hornjohnson.85}
R.~A. Horn and C.~R. Johnson.
\newblock \emph{Matrix Analysis}, Cambridge University Press, 1985.


\bibitem[4]{planque.18}
R.~Planqu\'e, J.~Hulshof, J.~C. Hendriks, B.~Teusink, and F.~J. Bruggeman.
\newblock Maintaining maximal metabolic rate using gene expression control.
\newblock \emph{PLoS Comp. Biol.}, 14\penalty0 (9):\penalty0 e1006412, 2018.

\bibitem[5]{planque.23}
R.~Planqu\'e, J.~Hulshof, and F.~J. Bruggeman.
\newblock How do bacteria control their metabolic state to maximise growth rate in varying environments?
\newblock \emph{forthcoming}, 2023.

\bibitem[6]{smillie.84}
  J. Smillie. \newblock Competitive and cooperative tridiagonal systems of differential equations. \newblock
  \emph{SIAM J. Math. Analysis},
  15\penalty0 (3):\penalty0 530-534. 1984.
  
\bibitem[7]{bruggeman.20}
F.~J. Bruggeman, R. Planqu\'e, D. Molenaar and B. Teusink. \newblock Searching for principles of microbial physiology. \newblock \emph{FEMS Microbiology Reviews}, 44\penalty0 (6):\penalty0  821-844. 2020.
 


\bibitem[8]{noor.16}
E. Noor, A. Flamholz, A. Bar-Even, D. Davidi, R. Milo and W. Liebermeister. 
\newblock The protein cost of metabolic fluxes: prediction from enzymatic rate laws and cost minimization. 
\newblock \emph{PLoS Comp. Biol.}, 12\penalty0 (11):\penalty0 e1005167. 2016.

\end{thebibliography}
\end{document}